\documentclass[12pt,leqno]{amsart}
\usepackage[dvips]{graphicx} \usepackage{amsfonts} \usepackage{amssymb} \usepackage{amsmath} \usepackage{amscd}
\usepackage{graphicx} 

\textheight9in \def\DATE{ 4 Chevat 5777 } \textwidth6.5in \hoffset-1.35cm \voffset-1cm 
\parindent10pt \footnotesep2mm \overfullrule3pt \newtheorem{theorem}{Theorem}
\newtheorem{definition}[theorem]{Definition} \newtheorem{corollary}[theorem]{Corollary}
 
 \newtheorem{lemma}[theorem]{Lemma}
 \newtheorem{proposition}[theorem]{Proposition}
  \newcommand\C{\mathbb{C}}
   
\newcommand\K{\mathbb{K}}  
     
\newcommand\ds{\displaystyle} 
\newcommand\pf{\noindent{\it Proof. }} 
\newcommand\ra{\rightarrow}
\email{elisabeth.remm@uha.fr}

\begin{document} 
\title{3-dimensional  skew-symmetric algebras and the variety of Hom-Lie algebras} 
\author{Elisabeth Remm}
  \date{ }
 \maketitle
 \begin{center}
Universit\'e de Haute Alsace, 4 rue des Fr\`eres Lumi\`ere, F68093 Muhouse.

\end{center}


\begin{abstract} An algebra is called skew-symmetric if its multiplication operation is a skew-symmetric bilinear application. We determine all these algebras in dimension $3$ over a field of characteristic different from $2$. As an application, we determine the subvariety of $3$-dimensional Hom-Lie algebras. For this type of algebras, we study also the dimension $4$. 
  \end{abstract}


\medskip


\bigskip

\noindent  INTRODUCTION. An algebra over a field $\K$ is a $\K$-vector space  equipped with a bilinear product. 
The multiplication operation in an algebra may or may not be associative, leading to the notions of associative algebras and nonassociative algebras. If we denote by $\mu$ this multiplication, we do not assume that $\mu$ satisfies some quadratic, ternary, or $n$-ary relations.  When we consider only the finite dimensional framework, that is when the $\K$-vector space is of finite dimension, one of the first natural problem which then arise is the determination of all these algebras. For example, the classification up to isomorphism of algebras for a given dimension  seems interesting and it looks at first easy to solve. It is strange that this problem was solved only in dimension $2$ (see the works of Goze-Remm, or Peterssen or Bekbaev (\cite{ Be, GRdim2, Pet}). For special classes of algebras, this classification work was carried  on for greater dimensions, but always relatively small. For example, associative algebras are classified up to dimension $6$, Lie algebras up to dimension $7$ (\cite{GozeKhakimdjanov}). Of course, the classification is complete for some special classes, for example the simple Lie algebras, the simple associative algebras. But for the general case, we have not many informations, probably because we know very few invariants in multilinear algebra.
Furthermore, when the classification exists (always with the fear that this classification could be incomplete), it is often difficult to use. Consider for example the classification of the complex nilpotent Lie algebras  of dimension $7$. We can think that this list is complete. There are in this case $5$  one-parameter families and more than one hundred of algebras (this number depends of the authors and it is not an invariant). Then it is often difficult to recognize among this list a given Lie algebra especially when if it is not written in a basis respecting the invariants used to obtain the classification. And it is not convenient to determine subclasses, for example to determine the contact nilpotent Lie algebras. To test each  algebra of this list is more than boring. That is why, as we presented it in previous works, we prefer to determine classes  invariant  by isomophism and minimal in a sense that we specify. 

We are interested in all this work in the $3$-dimensional algebras over a field $\mathbb{K}$ with skew-symmetric multiplication. We suppose in this work that $\K$ is an arbitrary field of characteristic different to $2$. The plane of this work is the following: first of all we study the automorphism group of a $3$-dimensional algebra. In particular, we characterize among all these algebras those which are Lie algebras  by studying the dimension of the automorphism group. In the second part, we classify the nilpotent 
skew-symmetric algebras by showing that all these algebras are Lie algebras (but it is not the case in greater dimension).
Next we classify the solvable skew-symmetric algebras and, to end this classification, the non-solvable case. The last section is essentially dedicated to the study of Hom-Lie algebras, which are a particular class of skew-symmetric algebras. We show that any skew-symmetric algebra of dimension $3$ is Hom-Lie. This is no longer true  in dimension $4$ where we study the algebraic variety of Hom-Lie algebras.

\medskip \section{The automorphism group of a $3$-dimensional skew-symmetric algebra}

\subsection{Generalities} An algebra over a field $\K$ is a $\K$-vector space $V$ with a multiplication operation given by a
bilinear map $$\mu: V \times V\rightarrow  V.$$ Such a multiplication is a linear tensor on $V$ of type
$(2,1)$, that is $2$-times contravariant and $1$-time covariant. 
In all this paper, the vector space $V$ is finite-dimensional and fixed. We can  for example consider $V=\K^n$. We denote by $A=(V,\mu)$ a $\K$-algebra with
multiplication $\mu$. 
We also assume in
this work that the field $\K$ is of characteristic  zero and that the multiplication $\mu$ is skew-symmetric,
that is $$\mu(X,Y)=-\mu(Y,X).$$ In this case, we shall say that $A$ is a skew-symmetric algebra, and we shorten the name by writing ss-algebra. 
Two $\K$-algebras $A=(V,\mu )$ and $A'=(V,\mu ')$ are  isomorphic if there is a linear
isomorphism $$f : V\rightarrow  V$$ such that $$f(\mu (X,Y))=\mu '(f(X),f(Y))$$ for all $X,Y \in V$ and we denote by $Aut(A)$ or $Aut(\mu)$ the group of automorphisms of the algebra $A=(V,\mu)$.

Assume now that $\dim V=3$.
Let $\{e_1,e_2,e_3\}$ be a fixed basis of $V$. A
general skew-symmetric bilinear map $\mu$ has the following expression 
\begin{equation}\label{Jac3} \left\{ \begin{array}{l}
\mu(e_1,e_2)=\alpha_1e_1+\beta _1e_2+ \gamma_1e_3, \\ \mu(e_1,e_3)=\alpha_2e_1+\beta _2e_2+\gamma_2e_3, \\
\mu(e_2,e_3)=\alpha_3e_1+\beta _3e_2+\gamma_3e_3, \\
 \end{array}
 \right.
\end{equation}
 and the set of these applications is a $9$-dimensional vector space $ss\mathcal{A}lg_3$ parametrized by the structure constants $\alpha_i,\beta_i,\gamma_i$, $i=1,2,3$, identifying $\mu$ with its structure constants. The linear group $GL(V)$ acts on
 $ss\mathcal{A}lg_3$:
 $$(g,\mu) \in GL(V) \times ss\mathcal{A}lg_3 \ra \mu_g \in
ss\mathcal{A}lg_3$$ with $$\mu_g(X,Y)=g^{-1}\mu (g(X),g(Y))$$ for any $X,Y \in V.$
Let $\mathcal{O}(\mu)$ be the orbit of $\mu$ corresponding to this action. It is a nonsingular algebraic subvariety of $ss\mathcal{A}lg_3$ (see for example \cite{Mar}). To compute the dimension of the subvariety $\mathcal{O}(\mu)$, it is sufficient to compute the dimension of its tangent space at the point $\mu$. Since $ss\mathcal{A}lg_3$ is linear, its tangent space $T_\mu(ss\mathcal{A}lg_3)$ at $\mu$ is identified to  
itself and it is constituted by skew-symmetric bilinear map $\psi \in ss\mathcal{A}lg_3$. The tangent space $T_\mu(\mathcal{O}(\mu))$ at the point $\mu$ to the orbit is a linear subspace of $T_\mu(ss\mathcal{A}lg_3)=ss\mathcal{A}lg_3$ whose elements are the skew-symmetric applications of type $\psi=\delta_\mu f$ with $f \in End(V)$ and
$$\delta_\mu
 f(X,Y)=\mu(f(X),Y)+\mu(X,f(Y))-f(\mu(X,Y))$$
 for any $X,Y \in V$.  Recall also that $f$ is a derivation of the algebra $A=(V,\mu)$ is $\delta_\mu(f)=0$. The set $Der(\mu)$ of derivations of $A$ is a Lie algebra which is the Lie algebra of the algebraic group $Aut(\mu)$.

 \subsection{The matrix $M_\mu$ associated to $\mu$}

 To compute the dimension of $T_\mu(\mathcal{O}(\mu))$, we shall consider the following  matrix
 $M_\mu$ associated to $\mu$, which will give a linear representation of $\mu$.  Recall that we have fixed a basis  $\{e_i\}$ a basis of $V$ and identified the skew-symmetric applications on $V$ with its structure constants related to this basis. Let $f \in gl(V)$ and
we denote also by $f$ the matrix of $f$ related to the fixed basis. If $f=(f_{ij})$, we
 consider the vector $v(f)$ of $\K^{9}$: $$v_f=(f_{1,1},f_{2,1},f_{3,1},f_{1,2},\cdots,f_{2,3},f_{3,3}).$$
We have
  $$\delta f(e_i,e_j)= \mu(f(e_i),e_j)+\mu(e_i,f(e_j))-f(\mu(e_i,e_j)$$
  for $1 \leq i < j \leq 3$ and we still denote $\delta f(e_i,e_j)$ the column matrix o its components in the basis $\{e_1,\cdots, e_3\}$.
 The vector  $$\left(
  \begin{array}{l}
     \delta f(e_1,e_2)   \\
     \delta f(e_1,e_3)   \\
\delta f(e_2,e_3)   \\ \end{array}
 \right)
 $$ is then a matrix in $\mathcal{M}_{9,1}.$ It 
 corresponds to a matrix product:
 $$M_\mu \cdot \, ^tv_f.$$
 This matrix $M_\mu$ is a square matrix of order $9$ and the map
 $$\mu \ra M_\mu$$
 gives a linear representation of $\mu$.

 If we consider the  structure constants of $\mu$ related to the basis $\{e_i\}$ given in
 (\ref{Jac3}),
then  $M_\mu$ is equal to
 $$M_\mu=\begin{pmatrix}
     0 & 0  & -\alpha_3 &-\beta_1 & \alpha_1  &\alpha_2 & -\gamma_1 &0 &  0  \\
    \beta_1  &-\alpha_1  & -\beta_3 & 0& 0&\beta_2 &0 &-\gamma_1 &0 \\
   \gamma_1 & 0 & -\gamma_3-\alpha_1 & 0 &\gamma_1 & \gamma_2-\beta_1 &0 & 0 & -\gamma_1\\
      0 & \alpha_3  & 0 &-\beta_2 & 0 & 0 & -\gamma_2 &\alpha_1 &  \alpha_2  \\
    \beta_2  &-\alpha_2+\beta_3  & 0& 0& -\beta_2 & 0  &0 &-\gamma_2+\beta_1 &\beta_2 \\
   \gamma_2 & \gamma_3 & -\alpha_2 & 0  & 0 & -\beta_2 &0 & \gamma_1 & 0\\
   - \alpha_3 & 0 & 0 &\alpha_2-\beta_3 & \alpha_3 & 0 & -\alpha_1-\gamma_3 & 0 & \alpha_3 \\
     0    & -\alpha_3& 0 & \beta_2 &0 &0 &-\beta_1 & -\gamma_3 &\beta_3 \\
 0 & 0 & -\alpha_3 & \gamma_2 & \gamma_3 & -\beta_3 & -\gamma_1 & 0 & 0 \\
\end{pmatrix}
 $$
 \begin{proposition}
 For any $\mu$, the matrix $M_\mu$ is singular. Moreover, there exists $\mu$ such that ${\rm Rank}(M_\mu)=8.$
 \end{proposition}
 \pf The first part is a direct computation of the determinant. We can see also that $0$ is always an eigenvalue.
 Now we consider $\mu$ given by $(\beta_1=1,\gamma_2=2,\alpha_3=1)$. For this algebra the kernel of $M_\mu$ is of
 dimension $1$ and constituted of matrices
 $$f=
 \begin{pmatrix}
     0 & 0 & 0   \\
    0  & b_2 & 0 \\
    0 & 0 & -b_2
\end{pmatrix}
 $$
and ${\rm Rank}(M_\mu)=8.$ For a generic algebra, that is without algebraic conditions between the structure
constants, the kernel of $M_\mu$ is a $1$-dimensional vector space generated by the vector

$$ \begin{array}{l}
    ( -\alpha_1 \alpha_2 \beta_1 + \alpha_1^2 \beta_2 - \alpha_2^2 \gamma_1 + \beta_3^2 \gamma_1 + \alpha_1
    \alpha_2 \gamma_2 - \beta_1 \beta_3 \gamma_3 +
    \beta_3 \gamma_2 \gamma_3 - \beta_2 \gamma_3^2,\\
   \alpha_2 \beta_1^2 + \alpha_1 \beta_1 \beta_2 - 2  \beta_2 \gamma_1(\alpha_2+ \beta_3 ) +\alpha_2 \beta_1
   \gamma_2 +\alpha_1 \beta_2 \gamma_2 + \beta_1 \beta_3 \gamma_2 - \beta_3 \gamma_2^2 +\beta_1 \beta_2 \gamma_3
   +\beta_2 \gamma_2 \gamma_3,\\
   2\beta_2 \gamma_1 (\alpha_1 -  \gamma_3 ) - \alpha_2 \beta_1 \gamma_1 + \beta_1 \beta_3 \gamma_1 - \alpha_1
   \beta_1 \gamma_2 - \alpha_2 \gamma_1 \gamma_2 + \beta_3 \gamma_1 \gamma_2 +
     \alpha_1 \gamma_2^2 - \beta_1^2 \gamma_3 + \beta_1 \gamma_2 \gamma_3,  \\
       \alpha_1 \alpha_3 \beta_1 - \alpha_1^2 \beta_3 + 2\alpha_3 \gamma_1 (\alpha_2  +  \beta_3 ) - \alpha_1
       \alpha_3 \gamma_2 -
    \alpha_1 \alpha_2 \gamma_3 - \alpha_3 \beta_1 \gamma_3 - \alpha_1 \beta_3 \gamma_3 + \alpha_3 \gamma_2
    \gamma_3 - \alpha_2 \gamma_3^2,\\
    \alpha_3 \beta_1^2 - \alpha_1 \beta_1 \beta_3 - \alpha_2^2 \gamma_1 + \beta_3^2 \gamma_1 +
     \alpha_1 \alpha_2 \gamma_2 - \alpha_3 \gamma_2^2 - \beta_1 \beta_3 \gamma_3 + \alpha_2 \gamma_2 \gamma_3,
     \\ \alpha_1 \beta_3 \gamma_1- \alpha_1 \alpha_2 \gamma_1 - 2 \alpha_3 \gamma_1( \beta_1  + \gamma_2) +
     \alpha_1^2 \gamma_2  +
    \alpha_1 \beta_1 \gamma_3 + \alpha_2 \gamma_1 \gamma_3 - \beta_3 \gamma_1 \gamma_3 + \alpha_1 \gamma_2
    \gamma_3 + \beta_1 \gamma_3^2, \\
    \alpha_2 \alpha_3 (\beta_1 -  \gamma_2 )- 2\alpha_3 \beta_2( \alpha_1 -  \gamma_3) + \alpha_1 \alpha_2
    \beta_3 +\alpha_3 \beta_1 \beta_3 -
       \alpha_1 \beta_3^2 - \alpha_3 \beta_3 \gamma_2 + \alpha_2 \gamma_3(\alpha_2   -\beta_3),\\
   \alpha_1 \alpha_2 \beta_2 -\alpha_2^2 \beta_1   + \alpha_2 \beta_1 \beta_3 +\alpha_1 \beta_2 \beta_3 -
    2 \alpha_3 \beta_2( \gamma_2 -  \beta_1)+ \alpha_2 \beta_3 \gamma_2 - \beta_3^2 \gamma_2 + \alpha_2 \beta_2
    \gamma_3 + \beta_2 \beta_3 \gamma_3, \\
    \alpha_1 \alpha_2 \beta_1 + \alpha_3 \beta_1^2 - \alpha_1^2 \beta_2 - \alpha_1 \beta_1 \beta_3 - \alpha_3
    \gamma_2^2 + \alpha_2 \gamma_2 \gamma_3 -
    \beta_3 \gamma_2 \gamma_3 + \beta_2 \gamma_3^2).
\end{array} $$

\begin{corollary}For any $3$-dimensional ss-algebra $A=(V,\mu)$  over $\K$, the  dimension of the automorphism group
$Aut(\mu)$ is greater than or equal  to $1$. 
\end{corollary} 
\pf In fact, $Aut(\mu)$ is an algebraic group whose
Lie algebra is isomorphic to $Der(\mu)$, the Lie algebra of derivations of $\mu$, that is the subspace of
$gl(V)$ whose element $f$ satisfy $\delta (f) = 0$ or equivalently the vector $^tv_f$ is in the kernel of
$M_\mu$. Since the rank of this matrix is smaller that $8$, its kernel is bigger than $1$.

\medskip

\begin{proposition} The algebra $A=(V,\mu)$ is a Lie algebra if and only if there exists a non zero vector $X$ in
$V$  such that the linear endomorphism $L_X: Y \ra \mu(X,Y)$ is in the kernel of $M_\mu$. 
\end{proposition} 
\pf If $f \in End(A)$ then, considering the basis $\{e_i\}$, we identify $f$ with the matrix $$
  \begin{pmatrix}
     a_1 &   b_1 & c_1 \\
        a_2 &   b_2 & c_2 \\
           a_3 &   b_3 & c_3 \\
\end{pmatrix} $$ and we associate it with  the vector $$v_f=\, ^t(a_1,a_2,a_3,b_1,b_2,b_3,c_1,c_2,c_3).$$
 If we take $X=e_1$, then
$$L_{e_1}= \begin{pmatrix}
     0 &  \alpha_1 & \alpha_2 \\
        0 &   \beta_1 & \beta_2 \\
           0&   \gamma_1 & \gamma_2 \\
\end{pmatrix} $$ and $M_\mu(v_{L_{e_1}})=0$ is equivalent to $$ \left\{ \begin{array}{l}
   -\alpha_1\beta_3 + \alpha_3\beta_1  -\gamma_3\alpha_2 + \alpha_3\gamma_2=0, \\
     \beta_2\alpha_1 -\beta_1\alpha_2  -\gamma_3\beta_2 +\beta_3\gamma_2=0, \\
\gamma_2\alpha_1 + \gamma_3\beta_1  -\beta_3\gamma_1 -\gamma_1\alpha_2 = 0 \\ \end{array} \right. $$ which
corresponds to Jacobi's conditions.

\begin{proposition} If the ss-algebra $A=(V,\mu)$ is a Lie algebra, then the rank of  $M_\mu$ is smaller than or equal
to $6$. \end{proposition} 
\pf  In fact, since $A=(V,\mu)$ is a Lie algebra, the endomorphisms $L_X(Y)=\mu(X,Y)$ are
derivations and the vectors $v_{L_X}$ are in the kernel of $M_\mu$. If $L_{e_1},L_{e_2},L_{e_3}$ are linearly
independent, then $\ker M_\mu \geq 3$ and the rank of $M_\mu $ is smaller than $6$. If these vectors are linearly
dependent, we have a non trivial linear combination between these vectors and without loss of  generality, we can
assume that $L_{e_1}=0$. But $v_{L_{e_1}}=(0,0,0,\alpha_2,\beta_2,\gamma_2,\alpha_3,\beta_3,\gamma_3)$ and
$v_{L_{e_1}}=0$ implies that 
$$M_\mu=\begin{pmatrix}
     0 & 0  & 0 &-\beta_1 & \alpha_1  &0 & -\gamma_1 &0 &  0  \\
    \beta_1  &-\alpha_1  & 0 & 0& 0&0 &0 &-\gamma_1 &0 \\
   \gamma_1 & 0 & -\alpha_1 & 0 &\gamma_1 &-\beta_1 &0 & 0 & -\gamma_1\\
      0 &0  & 0 &0 & 0 & 0 &0 &\alpha_1 &0  \\
  0 &0 & 0& 0&0 & 0  &0 &\beta_1 &0\\
 0& 0 &0 & 0  & 0 & 0 &0 & \gamma_1 & 0\\
   0 & 0 & 0 &0 & 0 & 0 & -\alpha_1 & 0 &0\\
     0    & 0& 0 & 0 &0 &0 &-\beta_1 & 0 &0 \\
 0 & 0 & 0 & 0 & 0& 0 & -\gamma_1 & 0 & 0 \\
\end{pmatrix}
 $$
 which  rank is less or equal to  $3$.

\medskip

 \noindent The converse of this proposition is not true. We can find non Lie algebras whose associated matrix
 $M_\mu$ is of rank $6$. For example we consider the algebra
$$ \left\{ \begin{array}{l} \mu(e_1,e_2)=e_2, \\ \mu(e_1,e_3)=\gamma_2e_3, \\ \mu(e_2,e_3)=e_1. \\
 \end{array}
 \right.
$$ This algebra is a Lie algebra if and only if $\gamma_2=-1.$ In this case the rank is $6$ and the kernel is
composed of the endomorphisms $$\begin{pmatrix}
      0&-a_3& a_2    \\
     -a_2& b_2  & 0\\
    a_3 & 0 & -b_2
\end{pmatrix} $$ If $\gamma_2=1$, the algebra is not a Lie algebra but the rank is also $6$. In this case the
kernel is composed of endomorphisms $$\begin{pmatrix}
      0&0& 0    \\
     0& b_2  & c_2\\
 0 & b_3 & -b_2
\end{pmatrix} $$ While this last algebra is not a Lie algebra, we can embedded this algebra in a larger class of
algebras such as algebras of Lie type (see \cite{Ma}) or Hom-Lie algebras. We shall see that in the last section.

 \subsection{What happen in dimension $4$?}
 If $A=(V_4,\mu)$ is a $4$-dimensional skew-symmetric $\K$-algebra, then $M_\mu$ is a non square matrix of order $24 \times 16$. It
 is not too much complicated to write this matrix, but we do not make here. If we identified $A$ with $\mu$ and $\mu$ with its structure constant related with a fixed basis $\{e_1,e_2,e_3,e_4\}$ of $V$, the set of $4$-dimensional ss-algebras is a $24$-dimensional vector space  denoted $ss\mathcal{A}lg_4$.

 \begin{proposition}
 There exists a Zariski-open set in the affine variety $ss\mathcal{A}lg_4$ 
 constituted of ss-algebras $(V,\mu)$ whose automorphism group $Aut(\mu)$ is of dimension $0$.
  \end{proposition}
  \pf In fact, for a generic algebra of $ss\mathcal{A}lg_4$, that is without relations between the structure
  constants, the corresponding matrix $M_\mu$ is of rank $16$. Then, any derivation is trivial and the algebraic
  group $Aut(\mu)$ is of dimension $0$.
  For example, for the following algebra
$$ \left\{ \begin{array}{l} \mu(e_1,e_2)=e_1+e_4, \\ \mu(e_1,e_3)=e_2+e_3,\\ \mu(e_1,e_4)=e_3+e_4, \\
\mu(e_2,e_3)=e_1+e_3, \\ \mu(e_2,e_4)=e_2,\\ \mu(e_3,e_4)=e_2+e_3, \\
 \end{array}
 \right.
$$ satisfies this property.

Let us note that $ss\mathcal{A}lg_4$ is an affine algebraic variety isomorphic to $\K^{24}$ and the Zariski open
set  constituted of algebras $A=(V_4,\mu)$ whose automorphism group $Aut(\mu)$ is of dimension $0$ is a finite union
of algebraic components where each one is the complementary of an hypersurface, then they are algebraic
subvarieties of $ss\mathcal{A}lg_4$.

\medskip

\noindent \noindent{\bf{Remark: Case of Dimension $2$.}}
 Before studying the dimension $3$, we have naturally studied the  $2$-dimensional case  in \cite{GRKhukhro}. In
 this case, if $\mu(e_1,e_2)=\alpha e_1+\beta e_2$, then
 $$M_\mu=
 \begin{pmatrix}
    0  &  0 & -\beta &\alpha   \\
    \beta  &  -\alpha & 0 & 0
\end{pmatrix} $$ Its rank is $2$ as soon as $\alpha \neq 0$ or $\beta \neq 0$. In this case $\ker (M_\mu)$ is of
dimension $2$. This implies that the group of automorphisms of $A=(V,\mu)$ is an algebraic group of dimension $2$.
We have also noted that  any ss-algebra of dimension $2$ is a Lie algebra.

\subsection{Deformations and rigidity}
A deformation of $\mu$, in the Gerstenhaber'sense, is given by a tensor $(2,1)$ in
$\K[t] \otimes ss\mathcal{A}lg_3.$ Grosso-modo, the notion of deformation permits to describe a neighborhood of
$\mu$ in $ss\mathcal{A}lg_3$. There exists
a cohomological approach of the notion of deformations. In our situation, the complex of cohomology is well
known. The space of $p$-cochains is the space of $p$-linear skew symmetric applications, the first space of
cohomology corresponds to the Kernel of $M_\mu$, and the second space is a factor space isomorphic to $T_\mu(ss\mathcal{A}lg_3)/ Im
M_\mu$. This complex is well described in the context of operads. We denote by $\mathcal{S}ign$ the quadratic
operad encoding the category of skew-symmetric algebras. Recall its construction (see also (see \cite{MR}).
Let $\K[\Sigma_2]$ be the group algebra of the symmetric group of degree 2. Considered as a $\Sigma_2$-module we
have that $\K[\Sigma_2]=1\! \!1_2 \oplus sgn_2$ where  $1\! \!1$ is the one-dimensional representation and
$sgn_2$ the one-dimensional signum representation. We then consider $\Gamma(sgn_2)$ the free operad generated by
a skew symmetric operation. For this operad, we have in particular $\Gamma(sgn_2)(0)=0$,
$\Gamma(sgn_2)(1)\simeq \K,$  $\dim \Gamma(sgn_2)(2) =1$ because $ \Gamma(sgn_2)(2)$ is the $\K[\Sigma_2]$-module
generated by $\{ x_1x_2=- x_2x_1\}$,  $\dim \Gamma(sgn_2)(3)=3$, $ \Gamma(sgn_2)(3)$ is the $\K[\Sigma_2]$-module
generated by $\{ (x_1x_2) x_3), (x_2 x_3)x_1,(x_3x_1) x_2)\},$$\dim  \Gamma(sgn_2)(4)=15$ because $
\Gamma(sgn_2)(2)$ is the $\K[\Sigma_2]$-module generated by

$\begin{array}{l} \{(x_1x_2)(x_3x_4),(x_1x_3)(x_2x_4), (x_1x_4)(x_2x_3), ((x_1x_2)x_3)x_4, ((x_1x_2)x_4)x_3,
((x_1x_3)x_2)x_4, \\
((x_1x_3)x_4)x_2,((x_1x_4)x_2)x_3,((x_1x_4)x_3)x_2,((x_2x_3)x_1)x_4,((x_2x_3)x_4)x_1,((x_2x_4)x_1)x_3,\\
((x_2x_4)x_3)x_1,((x_3x_4)x_1)x_2,((x_3x_4)x_2)x_1\}. \end{array} $

\noindent This operad is a Koszul operad, then the cohomology which parametrizes the deformations is the operadic
cohomology. 

\medskip

Recall also the notion of rigidity which traduces the fact that any deformations of $\mu$ is isomorphic to $\mu$.
\begin{definition}
A skew-symmetric algebra $A=(V,\mu)$ is rigid if the orbit $\mathcal{O}(\mu)$ is Zariski-open in the algebraic variety $ss\mathcal{A}lg_3$.
\end{definition}
This topological notion can be repalced by an algebraic condition. If we denote by $H^*_S(\mu,\mu)$ the complex of deformations, since the variety $\mathcal{A}lg_3$ is an affine
reduced variety, an element $\mu \in \mathcal{A}lg_3$ is rigid if and only if $H^2_S(\mu,\mu)=0$. 
We deduce of the previous results

\begin{proposition} For any $\mu \in ss\mathcal{A}lg_3$, we have $$\dim \mathcal{O}(\mu) \leq 8.$$ In particular
none of  $3$-dimensional ss-algebra is rigid in $ss\mathcal{A}lg_3$. \end{proposition}
 \pf Since $ss\mathcal{A}lg_3$  is
a linear plane and then a reduced algebraic variety, and since $\dim ss\mathcal{A}lg_3- \dim  \mathcal{O}(\mu) \geq
1$, none of these algebras is rigid in $ss\mathcal{A}lg_3$.

This notion of rigidity which concerns an element can be extended to parametrized families of algebras.

\begin{definition} Let $F_{t_1, \cdots t_k} $ be a family of algebras of $ss\mathcal{A}lg_3$ parametrized by $t_1,
\cdots ,t_3.$ This family is rigid if its orbit by the action of the group $GL(V)$ is Zariski open in
$ss\mathcal{A}lg_3.$ 
\end{definition} 
If $\mathcal{O}(F_{t_1, \cdots t_k} )$ is this orbit, its rigidity implies
that $\overline{\mathcal{O}(F_{t_1, \cdots t_k} )}$ (the closure in the Zariski sense) is an algebraic component
of $ss\mathcal{A}lg_3.$ But $ss\mathcal{A}lg_3$ is connected so if $F_{t_1, \cdots t_k} $ is rigid,
$\overline{\mathcal{O}(F_{t_1, \cdots t_k} )}=\mathcal{A}lg_3.$

\begin{proposition} The family $F_{\beta_2,\gamma_2,\alpha_3,\beta_3,\gamma_3}$ whose elements are the algebras
(\ref{ns1}): $$ \left\{ \begin{array}{l} \mu(e_1,e_2)=e_3, \\ \mu(e_1,e_3)=\beta_2e_2+\gamma_2e_3, \\
\mu(e_2,e_3)=\alpha_3e_1+\beta_3e_2+\gamma_3e_3, \end{array} \right. $$ is rigid in $ss\mathcal{A}lg_3$.
\end{proposition}  
In fact any deformation of any element of this family is isomorphic to an element of this same family.

\section{Nilpotent case}

Let $(A,\mu)$ be a $\K$-algebra. We consider the descending central series $$\mathcal{C}^0(\mu)=A,
\mathcal{C}^1(\mu)=\mu(A,A), \  \mathcal{C}^k(\mu)=\mu( \mathcal{C}^{k-1}(\mu),A), \ k \geq 2.$$ The algebra $A$
is nilpotent if there exists $k$ such that $\mathcal{C}^k(\mu)=0.$ The smallest $k$ such that
$\mathcal{C}^k(\mu)=0$ is the nilindex. For a $3$-dimensional nilpotent algebra $(A,\mu)$, we have only the
following sequences: \begin{enumerate}
  \item $ A \supset \mathcal{C}^1(\mu)=0,$
  \item $ A \supset \mathcal{C}^1(\mu) \supset \mathcal{C}^2(\mu)=0,$
  \item $ A \supset \mathcal{C}^1(\mu) \supset \mathcal{C}^2(\mu)\supset \mathcal{C}^3(\mu)=0.$
\end{enumerate} It is obvious that the last term $\mathcal{C}^k(\mu)$ is contained in the center. Thus the first
case corresponds to abelian case, the last is impossible because  it would imply that $\dim
\mathcal{C}^1(\mu)=2$. Then it remains $ A \supset \mathcal{C}^1(\mu) \supset \mathcal{C}^2(\mu)=0.$ We consider
an adapted basis $\{e_1,e_2,e_3\}$. It satisfies \begin{equation} \label{Nil1} \left\{ \begin{array}{l }
  \mu(e_i,e_3)=0,    \\
       \mu(e_1,e_2)=\gamma_1 e_3,
\end{array} \right. \end{equation} with $\gamma_1 \neq 0.$ We deduce

\begin{proposition} Any $3$-dimensional nilpotent algebra is a Lie algebra. It is isomorphic to the abelian
3-dimensional Lie algebra or to the Heisenberg algebra $\mathfrak{h}_3$. \end{proposition}

\noindent Let us note that in the non abelian case the rank of $M_\mu$ is equal to $3$.

\medskip

\noindent {\bf Remark: What happen in dimension greater than or equal to $4$?} We have seen that in dimension
$3,$ any nilpotent algebra is a Lie algebra. This property doesn't extend to greater dimension. Let us consider
for example the filiform case, that is nilpotent algebra such that $$ A \supset \mathcal{C}^1(\mu) \supset
\cdots \supset \mathcal{C}^k(\mu)=0,$$ with $\dim A/ \mathcal{C}^1(\mu) =2$ and $\dim \mathcal{C}^k(\mu)/
\mathcal{C}^{k+1}(\mu) =1,\ i=1,\cdots,k-1.$ If in dimension less than or equal to $4$ these algebras are also
Lie algebras, in dimension 5 (or greater) we can find filiform non Lie algebras. Consider the family
\begin{equation} \left\{ \begin{array}{l}
\mu(e_1,e_i)=e_{i+1}, \ i=2,3,4,\\
     \mu(e_2,e_3)=ae_4+be_5, \\
     \mu(e_2,e_4)=ce_5 ,\\
      \mu(e_3,e_4)=de_5.
     \end{array}
\right. \end{equation}   Such an algebra is not a Lie algebra as soon as $a-c \neq 0$ or $d \neq 0$. But this algebra is provided with a Hom-Lie algebra structure (see Section 5).

\section{Solvable case}

Let $(A,\mu)$ a $\K$ algebra. We consider the derived series: $$\mathcal{D}^0(\mu)=A,
\mathcal{D}^1(\mu)=\mu(A,A), \  \mathcal{D}^k(\mu)=\mu( \mathcal{D}^{k-1}(\mu),\mathcal{D}^{k-1}(\mu)), \ k \geq
2.$$ The algebra $A$ is solvable if there exists $k$ such that $\mathcal{D}^k(\mu)=0.$ The smallest $k$ such that
$\mathcal{D}^k(\mu)=0$ is the solvindex.  For a 3-dimensional solvable algebra $(A,\mu),$ we have only the
following sequences: \begin{enumerate}
  \item $ A \supset \mathcal{D}^1(\mu)=0,$
  \item $ A \supset \mathcal{D}^1(\mu) \supset \mathcal{D}^2(\mu)=0,$
  \item $ A \supset \mathcal{D}^1(\mu) \supset \mathcal{D}^2(\mu)\supset \mathcal{D}^3(\mu)=0.$
\end{enumerate} In the first case $A$ is abelian. In the second case,  if $\dim A/ \mathcal{D}^1(\mu)=2$, we
consider an adapted basis $\{e_1,e_2,e_3\}$ and the hypothesis imply that \begin{equation} \label{SolLie1}
\left\{ \mu(e_i,e_j)=a_{ij}e_3,  \ 1 \leq i < j \leq 3. \right. \end{equation} and, not to consider again the
nilpotent case, we consider  $a_{13} \neq 0$ or $a_{23} \neq 0$. This algebra is a solvable Lie algebra for any
$a_{ij}$ and the rank
 of $M_\mu=5$. In particular the automorphism group is of dimension $4$.  Let us note that any algebras of
 (\ref{SolLie1}) is isomorphic to
$$ \left\{ \mu(e_1,e_3)=e_3, \right. $$ and in this case, the identity component of this group is constituted of
matrices
 $$\begin{pmatrix}
   1   &  0 & 0   \\
     x & y & 0 \\
     z & 0 & t
\end{pmatrix} $$ with $yt \neq 0.$

Assume now that $\dim A/ \mathcal{D}^1(\mu)=1$ and consider an adapted basis $\{e_1,e_2,e_3\}$. This means that
$\{e_2,e_3\}$ is a basis of $\mathcal{D}^1(\mu)$. By hypothesis $\mu(e_2,e_3)=0$. Then $\mu$ satisfies
\begin{equation} \label{SolLie2} \left\{ \begin{array}{l}
   \mu(e_1,e_2)=\beta_1e_2+\gamma_1e_3,\\
     \mu(e_1,e_3)=\beta_2e_2+\gamma_2e_3, \\
     \beta_1\gamma_2-\beta_2\gamma_1 \neq 0
\end{array} \right. \end{equation} Such algebra is also a Lie algebra. Recall that the classification, up to
isomorphism, of $3$-dimensional solvable Lie algebras is given in the book of Jacobson (\cite{Ja}). The rank of
$M_\mu$ is smaller or equal to $5$. For example, it is equal to $5$ for $$ \left\{ \begin{array}{l}
   \mu(e_1,e_2)=e_3,\\
     \mu(e_1,e_3)=e_2,\\
\end{array} \right. $$ and equal to $3$ for $$ \left\{ \begin{array}{l}
   \mu(e_1,e_2)=e_2,\\
     \mu(e_1,e_3)=e_3.\\
\end{array} \right. $$ 

\noindent Assume now that the descending sequence is filiform, that is $A \supset \mathcal{D}^1(\mu) \supset
\mathcal{D}^2(\mu)\supset \mathcal{D}^3(\mu)=0$. Let $\{e_1,e_2,e_3\}$ be an adapted basis to this flag.  We have
$$\mu(e_2,e_3)=\gamma_3e_3, \  \mu(e_1,e_2)=\beta_1e_2+\gamma_1e_3, \ \  \mu(e_1,e_3)=\beta_2e_2+\gamma_2e_3$$
with $ \gamma_3 \neq 0 $ and $\beta_1 \neq 0$ or $\beta_2 \neq 0$. Such algebra is never a Lie algebra. Let us
note that, since $\gamma_3 \neq 0$, giving $e'_2=\gamma_3^{-1}e_3$, we can consider that $\gamma_3=1$.

\begin{proposition} Any $3$-dimensional solvable algebra is a Lie algebra or  is isomorphic to the non-Lie
algebra given by \begin{equation} \label{sol} \left\{ \begin{array}{l}
     \mu(e_1,e_2)=\beta_1e_2+\gamma_1e_3,   \\
    \mu(e_1,e_3)=\beta_2e_2+\gamma_2e_3 \\
    \mu(e_2,e_3)=e_3,
\end{array} \right. \end{equation} with $\beta_1 \neq 0$ or $\beta_2\neq 0.$ \end{proposition} We have $${\rm
Rank}(M_\mu)=8$$ as soon as $\beta_2 \neq 0$ and the kernel is generated by
$$(-\beta_2,\beta_1+\gamma_2,-\beta_1^2-2\beta_2\gamma_1+\beta_1\gamma_2, 0,0,\beta_1,0,0,\beta_2).$$ If
$\beta_2=0$, then ${\rm Rank}(M_\mu)=7$ and the kernel is generated by $$\{(0, 0, -\gamma_1, 0, 0, 0, 0, 0, 1),
(0, 0, -\beta_1 + \gamma_2, 0, 0, 1, 0, 0, 0)\}.$$ It would be interesting to present this parameter $\beta_2$ as
an invariant. The following proposition gives an answer.

\begin{proposition} The Killing symmetric form $K_\mu(X,Y)=tr(L_X \circ L_Y)$ where $L_U(V)=\mu(U,V)$ of an
algebra (\ref{sol}) is nondegenerate if and only $\beta_2 \neq 0.$ \end{proposition} \pf The matrix of $K_\mu$ in
the basis $\{e_1,e_2,e_3\}$  is $$ \begin{pmatrix}
    \beta_1^2+2\beta_2\gamma_1+\gamma_2^2  &   \gamma_2 & -\beta_2 \\
      \gamma_2 &  1 & 0 \\
      -\beta_2 & 0 & 0
\end{pmatrix} $$ and its determinant is equal $-\beta_2^2.$ Let us note that $(A,\mu)$ is not a simple algebra,
$I=\K\{e_2,e_3\}$ is an ideal of $A$, but its Killing form is nondegenerate as soon as  $\beta_2 \neq 0.$ Let us
note also that $K_\mu$ is, in this case, a not invariant pseudo scalar product.

 \section{Non-Solvable case}
If $(A,\mu)$ is not solvable, there exists $k$ with $$ \mathcal{D}^k(\mu)=\mathcal{D}^{k-1}(\mu).$$ We have the
following possibilities: \begin{enumerate} \item $A=\mathcal{D}^1(\mu),$ \item $ A \supset
\mathcal{D}^1(\mu)=\mathcal{D}^2(\mu)\neq 0,$ \item $ A \supset \mathcal{D}^1(\mu) \supset
\mathcal{D}^2(\mu)=\mathcal{D}^3(\mu)\neq 0$. \end{enumerate} The sequence $ A \supset \mathcal{D}^1(\mu) \supset
\mathcal{D}^2(\mu)=\mathcal{D}^3(\mu)\neq 0$ is impossible because this implies that $\dim \mathcal{D}^2(\mu)=1$
and $\mu(\mathcal{D}^2(\mu),\mathcal{D}^2(\mu))\neq 0$, this is impossible since $\mu $ is skew-symmetric.
Similarly with the sequence $ A \supset \mathcal{D}^1(\mu)=\mathcal{D}^2(\mu)\neq 0.$ Thus, if $(A,\mu)$ is not
solvable, we have $$A=\mathcal{D}^1(\mu)$$ that is $$ \left\{ \begin{array}{l}
\mu(e_1,e_2)=\alpha_1e_1+\beta_1e_2+\gamma_1e_3, \\ \mu(e_1,e_3)=\alpha_2e_1+\beta_2e_2+\gamma_2e_3, \\
\mu(e_2,e_3)=\alpha_3e_1+\beta_3e_2+\gamma_3e_3 \end{array} \right. $$ with $$ \det\begin{pmatrix} \alpha_1
&\beta_1 &\gamma_1& \\ \alpha_2&\beta_2&\gamma_2 \\ \alpha_3&\beta_3&\gamma_3 \end{pmatrix} \neq 0. $$
\begin{lemma} There exists $X,Y$ independent in
$A$ such that $X,Y,\mu(X,Y)$ are also independent. \end{lemma} \pf Assume that this property is not true. In any
basis $\{e_1,e_2,e_3\}$ we must have $$ \left\{ \begin{array}{l} \mu(e_1,e_2)=\alpha_1e_1+\beta_1e_2, \\
\mu(e_1,e_3)=\alpha_2e_1+\gamma_2e_3, \\ \mu(e_2,e_3)=\beta_3e_2+\gamma_3e_3 \end{array} \right. $$ with $$
\det\begin{pmatrix} \alpha_1 &\beta_1 &0& \\ \alpha_2&0&\gamma_2 \\ 0 &\beta_3&\gamma_3 \end{pmatrix} =
-\alpha_1\beta_3\gamma_2-\alpha_2\beta_1\gamma_3\neq 0. $$ Let $e'_2=ae_2+be_3$. Then $rank
\{e_1,e'_2,\mu(e_1,e'_2)\}=2$ implies $\beta_1=\gamma_2.$ Likewise $rank \{e_2,e'_1=ae_1+be_3,\mu(e'_1,e_2)\}=2$
implies $-\gamma_3=\alpha_1$ and $rank \{e_3,e'_2=ae_1+be_2,\mu(e_3,e'_2)\}=2$ implies $\beta_3=\alpha_2.$ We
deduce $$ -\alpha_1\beta_3\gamma_2-\alpha_2\beta_1\gamma_3=-\alpha_1\alpha_2\beta_1+\alpha_2\beta_1\alpha_1=0.$$
This gives a contradiction.

From this lemma we can find a basis $\{e_1,e_2,e_3\}$ of $A$ such that $$ \left\{ \begin{array}{l}
\mu(e_1,e_2)=e_3, \\ \mu(e_1,e_3)=\alpha_2e_1+\beta_2e_2+\gamma_2e_3, \\
\mu(e_2,e_3)=\alpha_3e_1+\beta_3e_2+\gamma_3e_3, \end{array} \right. $$ and $\alpha_2\beta_3-\alpha_3\beta_2 \neq
0.$ We shall show that this family can be reduced to a family with only $5$ parameters. This condition implies
that $\alpha_2$ or $\alpha_3$ is non zero. We can assume that $\alpha_3 \neq 0$. In fact, if $\alpha_3=0$ and
$\alpha_2 \neq 0$ we have our result. If $\alpha_2 \neq 0$, then the change of basis $e'_1=e_1-\alpha_2/\alpha_3
e_2$ and $e'_i=e_i$ for $i=2,3$ permits to consider $\alpha_2=0$. We deduce:

\begin{proposition} Let $(A,\mu)$ a $3$-dimensional nonsolvable algebra. Then this algebra is isomorphic to the
algebra \begin{equation} \label{ns1} \left\{ \begin{array}{l} \mu(e_1,e_2)=e_3, \\
\mu(e_1,e_3)=\beta_2e_2+\gamma_2e_3, \\ \mu(e_2,e_3)=\alpha_3e_1+\beta_3e_2+\gamma_3e_3, \end{array} \right.
\end{equation} with $\alpha_3\beta_2 \neq 0,$ or \begin{equation} \label{ns2} \left\{ \begin{array}{l}
\mu(e_1,e_2)=e_3, \\ \mu(e_1,e_3)=\alpha_2e_1+\beta_2e_2+\gamma_2e_3, \\ \mu(e_2,e_3)=\beta_3e_2+\gamma_3e_3,
\end{array} \right. \end{equation} with $\alpha_2\beta_3 \neq 0$. \end{proposition} The algebras (\ref{ns1}) are
Lie algebras if and only if $\beta_3=\gamma_2=\gamma_3=0$, that is if we have $$ \left\{ \begin{array}{l}
\mu(e_1,e_2)=e_3, \\ \mu(e_1,e_3)=\beta_2e_2, \\ \mu(e_2,e_3)=\alpha_3e_1, \end{array} \right. $$ with
$\alpha_3\beta_2 \neq 0$ and for the algebras (\ref{ns2}) if and only if $\gamma_2=\gamma_3=0,
\beta_3=-\alpha_2$, that is if we have $$ \left\{ \begin{array}{l} \mu(e_1,e_2)=e_3, \\
\mu(e_1,e_3)=\alpha_2e_1+\beta_2e_2 , \\ \mu(e_2,e_3)=-\alpha_2e_2 \end{array} \right. $$ with $\alpha_2 \neq
0.$

For any generic algebra $\mu$ belonging to the family   (\ref{ns1}), the rank of $M_\mu$ is $8$, the kernel
being generated by the vector

$(-\beta_3^2 - \beta_3 \gamma_2 \gamma_3 + \beta_2 \gamma_3^2, 2 \beta_2 \beta_3 + \beta_3 \gamma_2^2 - \beta_2
\gamma_2 \gamma_3, -\beta_3 \gamma_2 + 2 \beta_2 \gamma_3, -
   2 \alpha_3 \beta_3 - \alpha_3 \gamma_2 \gamma_3, \beta_3^2 -
    \alpha_3 \gamma_2^2, $
    $ 2 \alpha_3 \gamma_2 + \beta_3 \gamma_3, \alpha_3 \beta_3 \gamma_2 - 2 \alpha_3 \beta_2 \gamma_3, 2 \alpha_3
    \beta_2 \gamma_2 + \beta_3^2 \gamma_2 -
    \beta_2 \beta_3 \gamma_3,
   \alpha_3 \gamma_2^2 + \beta_3 \gamma_2 \gamma_3 - \beta_2 \gamma_3^2)$

\noindent and of rank $6$ is case of Lie algebra, that is if
     $(\gamma_2,\beta_3,\gamma_3) =(0,0,0).$

\noindent We have a similar result for the family (\ref{ns2}), but in this case the kernel is generated by

$(\beta_3 \gamma_2\gamma_3 - \beta_2 \gamma_3^2,
  - \beta_3 \gamma_2^2 + \beta_2 \gamma_2 \gamma_3, 0, \alpha_2 \gamma_3^2
 , -\alpha_2 \gamma_2\gamma_3,
  0, \alpha_2^2\gamma_3 - \alpha_2 \beta_3\gamma_3
, \alpha_2 \beta_3 \gamma_2 - \beta_3^2 \gamma_2 + \alpha_2 \beta_2 \gamma_3 +
    \beta_2 \beta_3 \gamma_3, \alpha_2 \gamma_2\gamma_3 + \beta_3 \gamma_2\gamma_3 - \beta_2 \gamma_3^2.)$

\begin{proposition} Any $3$-dimensional nonsolvable algebra  is a simple algebra and belongs to the family
(\ref{ns1}) or (\ref{ns2}). For all these algebras, the automorphism group is of dimension $1$ except if the
algebra is a Lie algebra, in this case it is of dimension $3$\end{proposition}

 \section{Applications: Hom-Lie algebras}

\subsection{Hom-Lie algebras \cite{sil}} 

The notion of Hom-Lie algebras was introduced by Hartwig, Larsson and  Silvestrov  in \cite{sil}. Their principal
motivation concerns deformation of the Witt algebra. This Lie algebra is  the complexification of the Lie algebra
of polynomial vector fields on a circle. A basis for the Witt algebra is given by the vector fields
$$L_n=-z^{n+1}\ds \frac{\partial}{\partial z}$$
 for any  $n \in \mathbb {Z}.$
The Lie bracket  is given by $$[L_{m},L_{n}]=(m-n)L_{{m+n}}.$$ The Witt algebra is also viewed as  the Lie
algebra of derivations of the ring $\C[z,z^{-1}]$. Recall that a derivation on an algebra with product denoted by
$ab$ is a linear operator satisfying $D(ab)=D(a)b+aD(b)$. The Lie bracket of two derivations $D$ and $D'$ is
$[D,D']=D \circ D' - D' \circ D.$ We can also define a new class of linear operators generalizing derivations,
the Jackson derivate, given by $$D_q(f)(z)=\ds \frac{f(qz)-f(z)}{qz-z}.$$ It is clear that $D_q$ is a linear
operator, but its behavior on the product is quite different as the classical derivative:
$$D_q(fg(z))=g(z)D_q(f(z))+f(qz)D_q(g(z)).$$ The authors interpret this relation by putting \begin{equation}
\label{sigma} D_q(fg)=gD_q(f)+\sigma (f)D_q(g) \end{equation} where $\sigma$  is given by $\sigma(f)(z)=f(qz)$
for any $f \in \C[z,z^{-1}]$. Starting from (\ref{sigma}) and for a given $\sigma$, one defines a new space of
"derivations" on  $C[z,z^{-1}]$ constituted of linear operator $D$ satisfying this relation. With the classical
bracket we obtain is a new type of algebra so called $\sigma$-deformations of the Witt algebra.  This new
approach  leads  naturally to considerer the space of $\sigma$-derivations, that is, linear operators satisfying
(\ref{sigma}), to provide it with the multiplication associated with a bracket. This new algebra is not a Lie
algebra because the bracket doesn't satisfies the Jacobi conditions. The authors shows that this bracket
satisfies a "generalized Jacobi condition". They have called this new class of algebras the class of Hom-Lie
algebras. This notion, introduced in \cite{sil},
 since made the object of numerous studies and was also generalized (see \cite{MDL}).
 
We denote by $\mathcal{H}Lie_n$ the subset $ss\mathcal{A}lg_n$ whose elements are $n$-dimensional Hom-Lie algebras. We have seen that $ss\mathcal{A}lg_n$ is an affine variety isomorphic to $K^{n^2(n-1)/2}$. We shall study in particular the case $n=3$ and $n=4$, proving that in dimension $3$ any
skew-symmetric algebra is a Hom-Lie algebra and in dimension $4$,  $\mathcal{H}Lie_4$ is an algebraic
hypersurface in $ss\mathcal{A}lg_4$. We end this work by the determination of binary quadratic operads whose
associated algebras are Hom-Lie.

\begin{definition} A Hom-Lie algebra structure on the vector space $V$ is a triple $A=(V, \mu, \alpha)$ consisting of a
skew-bilinear map $\mu: V\times V \rightarrow V$ and
 a linear space homomorphism $f: V \rightarrow V$
 satisfying the Hom-Jacobi identity
$$
 \circlearrowleft_{x,y,z}{\mu(\mu(x,y),f(z))}=0
$$ for all $x, y, z$ in $A$, where $\circlearrowleft_{x,y,z}$ denotes summation over the cyclic permutations on
$x,y,z$. \end{definition}

For example, a Hom-Lie algebra whose endomorphism $f$ is the identity is a Lie algebra. We deduce, since any
$2$-dimensional skew-symmetric algebra (the multiplication $\mu$ is a skew-symmetric bilinear map)  is a Lie
algebra, that any $2$-dimensional algebra is a Hom-Lie algebra.

In the following section we are interested by the determination of all Hom-Lie algebras for small dimensions.

\subsection{Hom-Lie algebras of dimension $3$}

 Any  $3$-dimensional skew-symmetric $\K$-algebra $A=(V,\mu)$
is defined by its  structure constants $\{\alpha_i,\beta_i, \gamma_i \}_{i=1,2,3}$ with respect to a  given basis
$\{e_1,e_2,e_3\}:$ 
$$ \left\{ \begin{array}{l} \mu(e_1,e_2)=\alpha_1e_1+\beta _1e_2+ \gamma_1e_3, \\
\mu(e_1,e_3)=\alpha_2e_1+\beta _2e_2+\gamma_2e_3, \\ \mu(e_2,e_3)=\alpha_3e_1+\beta _3e_2+\gamma_3e_3. \\
 \end{array}
 \right.
 $$
Let $f$ be an element of  $gl(3,\K)$ and  consider its matrix in the same basis $\{e_1,e_2,e_3\}$ $$\left(
\begin{array}{lll}
   a_1 & b_1 & c_1   \\
     a_2 & b_2 & c_2   \\
   a_3 & b_3 & c_3   \\
\end{array} \right). $$ We then define the vector $$v_f=\ ^t (a_1,a_2,a_3,b_1,b_2,b_3,c_1,c_2,c_3).$$ 
For such
an algebra we associate the following matrix, $HL_{\mu}$,  belonging to $\mathcal{M}(3,9)$
the space of matrices of order $(3 \times 9)$ and given by 
 $$\begin{pmatrix}
   \alpha_1v_1+\alpha_2v_2+\alpha_3v_3\\
     \beta_1v_1+\beta_2v_2+\beta_3v_3 \\
     \gamma_1v_1+\gamma_2v_2+\gamma_3v_3
\end{pmatrix}
 $$
 where  $$\begin{array}{l}
   v_1=(-\beta_3,\alpha_3,0,\beta_2,-\alpha_2,0,-\beta_1,\alpha_1,0) ,   \\
    v_2=(-\gamma_3,0,\alpha_3,\gamma_2,0,-\alpha_2,-\gamma_1,0,\alpha_1),\\
    v_3=(    0,-\gamma_3,\beta_3,0,\gamma_2,-\beta_2,0,-\gamma_1,\beta_1).
\end{array} $$
Using the notation $ij.k$ in place of $\mu(\mu(e_i,e_i),e_k)$, we have 
$$HL_\mu=\begin{pmatrix}
23.1 & 23.2  & 23.3 & 31.1 & 31.2 & 31.3 & 12.1 & 12.2 & 12.3\\
  \end{pmatrix}
$$ with $$23.1=\begin{pmatrix}
    -\alpha_1 \beta_3 -\alpha_2\gamma_3  \\
   -\beta_1 \beta_3 -\beta_2\gamma_3  \\
    -\gamma_1 \beta_3 -\gamma_2\gamma_3  \\
\end{pmatrix} , \  23.2=\begin{pmatrix}
    \alpha_1 \alpha_3 -\alpha_3\gamma_3  \\
   \beta_1\alpha_3 -\beta_3\gamma_3  \\
    \gamma_1 \alpha_3 -\gamma_3\gamma_3  \\
\end{pmatrix} , \  23.3=\begin{pmatrix}
    \alpha_2 \alpha_3 +\alpha_3\beta_3  \\
   \beta_2\alpha_3 +\beta_3\beta_3  \\
    \gamma_2 \alpha_3 +\gamma_3\beta_3  \\
    \end{pmatrix}
$$ $$31.1=\begin{pmatrix}
    \alpha_1 \beta_2 +\alpha_2\gamma_2  \\
   \beta_1 \beta_2 +\beta_2\gamma_2 \\
    \gamma_1 \beta_2 +\gamma_2\gamma_2  \\
\end{pmatrix} , \  31.2=\begin{pmatrix}
    -\alpha_1 \alpha_2 +\alpha_3\gamma_2  \\
   -\beta_1\alpha_2 +\beta_3\gamma_2  \\
    -\gamma_1 \alpha_2 +\gamma_3\gamma_2  \\
\end{pmatrix} , \  31.3=\begin{pmatrix}
    -\alpha_2 \alpha_2 -\alpha_3\beta_2  \\
   -\beta_2\alpha_2 -\beta_3\beta_2  \\
    -\gamma_2 \alpha_2 -\gamma_3\beta_2  \\
    \end{pmatrix}
$$ $$12.1=\begin{pmatrix}
    -\alpha_1 \beta_1 -\alpha_2\gamma_1  \\
   -\beta_1 \beta_1 -\beta_2\gamma_1  \\
    -\gamma_1 \beta_1 -\gamma_2\gamma_1  \\
\end{pmatrix} , \  12.2=\begin{pmatrix}
    \alpha_1 \alpha_1 -\alpha_3\gamma_1  \\
   \beta_1\alpha_1 -\beta_3\gamma_1  \\
    \gamma_1 \alpha_1 -\gamma_3\gamma_1  \\
\end{pmatrix} , \  12.3=\begin{pmatrix}
    \alpha_2 \alpha_1 +\alpha_3\beta_1  \\
   \beta_2\alpha_1 +\beta_3\beta_1  \\
    \gamma_2 \alpha_1 +\gamma_3\beta_1  \\
    \end{pmatrix} .
$$

\begin{theorem} Any skew-symmetric $3$-dimensional algebra is a Hom-Lie algebra. \end{theorem} 
\pf Consider a basis
$\{e_1,e_2,e_3\}$ of the algebra $(A,\mu)$ and let $(\alpha_i,\beta_i,\gamma_i)$, $i=1,2,3$ be its structure
constants defined previously. Let $f$ be in $gl(3,\K)$ and let us consider the associated vector
$$v_f=(a_1,a_2,a_3,b_1,b_2,b_3,c_1,c_2,c_3).$$
 The endomorphism  $f\in gl(3,\K) $
 satisfies the Hom-Jacobi condition if and only if its corresponding vector $v_f=\ ^t
 (a_1,a_2,a_3,b_1,b_2,b_3,c_1,c_2,c_3)$ is in the kernel of the matrix $HL_\mu.$
 But this matrix   belongs to $\mathcal{M}(3,9)$ and represents a linear morphism
 $$t: \K^9 \rightarrow \K^3.$$
 From the rank theorem we have
 $$\dim {\rm Ker}\  t = 9- \dim {\rm Im}\  t \geq 6.$$
 Then this kernel is always non trivial and for any algebra $\mu$, there exists a non trivial element in the
 kernel. Then this algebra always admits a non trivial Hom-Lie structure.

 \medskip

\subsection{Classification of Hom-Lie algebras of dimension $3$} In the previous section we have
determinate  the $3$-dimensional skew-symmetric $\K$-algebra. Since any Hom-Lie algebra is
 skew-symmetric, we deduce the classification of Hom-Lie algebras. Moreover, for a given skew-symmetric algebra, we can calculate the endomorphisms $f$ associated with the Hom-Lie-Jacobi condition solving the linear system:
$$\begin{pmatrix} 23.1 & 23.2  & 23.3 & 31.1 & 31.2 & 31.3 & 12.1 & 12.2 & 12.3\\
  \end{pmatrix} ^tv_f=0.$$
For example,  the identity map whose associated vector is $$v_{Id}=(1,0,0,0,1,0,0,0,1)$$ is in
the kernel of $HL_\mu$ if and only if $\mu$ satisfies $$(23)1+(31)2+(12)3=0$$ that is if it is a Lie algebra.
We deduce
\begin{proposition}
If the Hom-Lie algebras $A=(V,\mu)$ and $A'=(V,\mu')$ are isomorphic, then 
the kernels of the associated matrices $HL_\mu$ and $HL_{\mu'}$ are isomorphic.
\end{proposition}

\subsection{Dimension $4$}

  Let $A=(V,\mu)$ a $4$-dimensional skew-symmetric $\K$-algebra.  Let us choose a basis $\{e_1,e_2,e_3,e_4\}$ of
  $V$ and let us consider the corresponding constants structure of $\mu$:
$$\mu(e_i,e_j)=\ds \sum_{k=1}^4 C_{i,j}^ke_k.$$
This algebra is a Hom-Lie algebra if there exists a linear endomorphism $f$ satisfying
the Hom-Lie Jacobi equations. The endomorphism $f$ is represented in the basis $\{e_1,e_2,e_3,e_4\}$ by a square
matrix $(a_{i,j}$ of order $4$. As in dimension $3$, we consider the vecteur $v_f \in \K^{16}$
  $$v_f=(a_{1,1},a_{2,1},a_{3,1},a_{4,1},a_{1,2},a_{2,2},\cdots,a_{3,4},a_{4,4}).$$
 Then $f$ satisfies the Hom-Lie conditions if and only if $v_f$ is solution on the linear system
 $$HL_\mu \cdot  ^tv_f=0$$
 where $HL_\mu$ is the square matrix of order $16$:
 $$
 HL_\mu=
 \begin{pmatrix}
   23.i  & 31.i & 12.i & 0    \\
     24.i &  41.i & 0 & 12.i \\
     34.i & 0 & 41.i & 13.i \\
     0 & 34.i & 42.i & 23.i
\end{pmatrix}
, \ \ i=1,2,3,4$$
and $ij.k$ is the matrix
$$
ij.k=\begin{pmatrix}
    \sum_l C_{i,j}^lC_{l,k}^1  \\
         \sum_l C_{i,j}^lC_{l,k}^2  \\
             \sum_l C_{i,j}^lC_{l,k}^3  \\
                 \sum_l C_{i,j}^lC_{l,k}^4  \\
\end{pmatrix}
$$
 The kernel of $HL_\mu$ is not trivial if and only if the rank of $HL_\mu$ is smaller or equal to $15$. 
  We deduce that $A=(V,\mu)$ is provided with a Hom-Lie structure if the structure constant $C_{i,j}^k$ satisfy the homogeneous polynomial equation of degree $16$:
$$\det (HL_\mu)=0.$$

  \begin{proposition}
  The set $\mathcal{HL}_4$ of $4$-dimensional $\K$-Hom-Lie algebras is provided with a structure of 
  algebraic variety embedded in $\K^{24}$.
  \end{proposition}

For any $A=(V,\mu) \in \mathcal{HL}_4$  we consider the vector space $\ker HL_\mu$. We thus define
a singular vector bundle $K(\mathcal{HL}_4)$ whose fiber over  $A=(V,\mu) $ is $\ker HL_\mu$. This
fiber corresponds to the set of Hom-Lie structure which can be defined on a given $4$-dimensional algebra.

 \noindent{\bf Remarks.} 
 \begin{enumerate}
  \item  In dimension $3$, $\mathcal{HL}_3$ is the affine variety $ss\mathcal{A}lg_3$ which is
 isomorphic to the affine space $\K^9$.  In this case, the fibers of $K(\mathcal{HL}_3)$ are
 vector spaces of dimension greater or equal to $6$.
 \item In dimension $4$, we are confronted with the resolution of the equation of degree $16$ with $24$ variables $\det (HL_\mu)=0.$ If $\K$ is algebraically closed, we can simplify this problem because the endomorphism $f$ admits a general reduced form
 $$\begin{pmatrix}
     a_{1,1} & a_{1,2} & 0 & 0    \\
      0& a_{2,2} & a_{2,3}& 0\\
      0 & 0 & a_{3,3} & a_{3,4} \\
      0 & 0 & 0 & a_{4,4}
\end{pmatrix} $$
In this case,  we have  to consider only  the reduced matrix $HL'_\mu$ of order $16 \times 7$:
$$\left( \begin{array}{cccccccccccccccc} 23 .1 & 31. 1 & 31.2 & 12. 2 & 12. 3  & 0 & 0\\ 24. 1 & 41. 1 & 41. 2 &
0 & 0 &   12. 3   & 12. 4\\ 34. 1 & 0 & 0 & 41. 2 & 41. 3 &  13. 3   & 13. 4\\ 0 & 34. 1 & 34. 2  &    42. 2 &
42. 3 & 23. 3   & 23. 4\\
  \end{array}
  \right)
$$
\end{enumerate}
It remains to verify that there exist at least one point of $ss\mathcal{A}lg_4$ which not belong to  $\mathcal{HL}_4$. 
Let us consider, for example, the
following algebra: $$
  \left\{
  \begin{array}{l}
     \mu(e_1,e_2)= e_2+2e_3-e_4   \\
     \mu(e_1,e_3	)= e_1+2e_2-e_3   \\
     \mu(e_1,e_4)= 2e_1-e_2+e_4   \\
     \mu(e_2,e_3	)= -e_1+e_3+2e_4   \\
     \mu(e_2,e_4	)= e_1+2e_2-e_3+3e_4   \\
     \mu(e_3,e_4)= -2e_1-e_2+e_3+2e_4   \\
\end{array} \right. $$ The associated matrix is $$ \left( \begin{array}{cccccccccccccccc}
 -5 & -1 & 3 & -4 & -1 & 1 & 1 & -6 & 0 & 3 & -3 & -3 & 0 & 0 & 0 & 0 \\
 0 & -3 & 0 & 0 & 0 & -1 & -2 & -4 & 6 & 2 & -1 & 0 & 0 & 0 & 0 & 0 \\
 1 & 3 & -1 & 1 & 5 & -3 & -1 & 3 & 0 & -3 & 2 & 1 & 0 & 0 & 0 & 0 \\
 -2 & -9 & -4 & 1 & -2 & -1 & -4 & -5 & -2 & -1 & 4 & 7 & 0 & 0 & 0 & 0 \\
 -5 & -4 & 4 & 6 & 2 & 1 & -5 & -3 & 0 & 0 & 0 & 0 & 0 & 3 & -3 & -3 \\
 3 & 3 & 4 & 4 & -2 & 0 & -5 & 4 & 0 & 0 & 0 & 0 & 6 & 2 & -1 & 0 \\
 -5 & 6 & 1 & -3 & -2 & -5 & 4 & -1 & 0 & 0 & 0 & 0 & 0 & -3 & 2 & 1 \\
 -1 & -8 & -3 & 5 & 2 & 5 & 4 & 1 & 0 & 0 & 0 & 0 & -2 & -1 & 4 & 7 \\
 -5 & -1 & 3 & -7 & 0 & 0 & 0 & 0 & 2 & 1 & -5 & -3 & 1 & -1 & -1 & 6 \\
 1 & -6 & -2 & -1 & 0 & 0 & 0 & 0 & -2 & 0 & -5 & 4 & 0 & 1 & 2 & 4 \\
 3 & -3 & -1 & 2 & 0 & 0 & 0 & 0 & -2 & -5 & 4 & -1 & -5 & 3 & 1 & -3 \\
 -3 & -6 & -6 & -3 & 0 & 0 & 0 & 0 & 2 & 5 & 4 & 1 & 2 & 1 & 4 & 5 \\
 0 & 0 & 0 & 0 & -5 & -1 & 3 & -7 & 5 & 4 & -4 & -6 & -5 & -1 & 3 & -4 \\
 0 & 0 & 0 & 0 & 1 & -6 & -2 & -1 & -3 & 5 & -4 & -4 & 0 & -3 & 0 & 0 \\
 0 & 0 & 0 & 0 & 3 & -3 & -1 & 2 & 5 & -6 & -1 & 3 & 1 & 3 & -1 & 1 \\
 0 & 0 & 0 & 0 & -3 & -6 & -6 & -3 & 1 & 8 & 3 & -5 & -2 & -9 & -4 & 1 \\
\end{array} \right) $$ and its determinant is not zero. This algebra is not Hom-Lie. As consequence we deduce that there exists an open set in $ss\mathcal{A}lg_4$ whose elements are not Hom-Lie algebras.

\begin{proposition} In the affine plane $ss\mathcal{A}lg_4$, there is a Zariski open set whose elements are
$4$-dimensional algebras without Hom-Lie structure. \end{proposition}

 \section{Algebras of Lie type of dimension $3$}

Since the multiplication $\mu$ is skew-symmetric, any quadratic relation concerning $\mu$ can be reduced to a
relation of type $$h_1(X,Y,Z)\mu(\mu(X,Y),Z)+h_2(X,Y,Z)\mu(\mu(Y,Z),X)+h_3(X,Y,Z)\mu(\mu(Z,X),Y)=0$$ where
$h_1,h_2,h_3$ are functions $h_i:V^3 \ra V$. In particular, we have the class of  algebras of Lie type
(\cite{Ma}) given by the relation $$\mu(\mu(X,Y),Z)+a(X,Y,Z)\mu(\mu(Y,Z),X)+b(X,Y,Z)\mu(\mu(Z,X),Y)=0$$ with $a
\neq 0$. For example, In the first section, we have considered the non Lie algebra

$$ \left\{ \begin{array}{l} \mu(e_1,e_2)=e_2, \\ \mu(e_1,e_3)=e_3, \\ \mu(e_2,e_3)=e_1. \\
 \end{array}
 \right.
$$ This algebra satisfies the identity $$\mu(\mu(X,Y),Z)+\mu(\mu(Y,Z),X)+b(X,Y,Z)\mu(\mu(Z,X),Y)=0$$ where $b$ is
a trilinear form satisfying $b(e_1,e_2,e_3)=-b(e_3,e_1,e_2)=1.$ It is an algebra of Lie type.  An important study
of the automorphism group of these algebras is presented in (\cite{Ma}). In this example, $\dim Aut(\mu)=3$ and
the identity component is the group of matrices $$\begin{pmatrix}
     1&0& 0    \\
     0& a  & b\\
 0 & c & d
\end{pmatrix} $$
 with $ad-bc=1$. It is isomorphic to $SO(2)$.
For all the non Lie algebras determined in the above sections, let us determine the  algebras of Lie type.
\begin{enumerate}
  \item Algebras (\ref{sol}).  Such algebras are of Lie type if $\beta_2=\beta_3=0$, that is if it is also a
      Lie algebra.
   \item Algebras (\ref{ns1}). Such algebras are never algebras of Lie type.
  \item  Algebras (\ref{ns2}). Such an algebra is a Lie   algebra  and only if it is a Lie algebra
      ($\gamma_2,\beta_3,\gamma_3)=(0,0,0))$.Assume $\alpha_2=1,\beta_2=0$.
We have a Lie algebra structure if and only if $(\gamma_2,\beta_3,\gamma_3)=(0,-1,0). $ We have a complex
structure of algebra of Lie type, which is not a Lie algebra if and only if $\gamma_2=\gamma_3=0$ and
$\beta_3$ is a complex root of $-1$ (we assume here that $\K=\C)$ that is $\beta_3=-j$ or $-j^2$. In this
case the relation is $$\mu(\mu(X,Y),Z)+j\mu(\mu(Y,Z),X)+j^2\mu(\mu(Z,X),Y)=0$$ or
$$\mu(\mu(X,Y),Z)+j^2\mu(\mu(Y,Z),X)+j\mu(\mu(Z,X),Y)=0.$$ In this case $M_\mu$ is of rank $8$ and the
automorphism group is of dimension $1$ and generated bu the automorphism $$\begin{pmatrix}
    u  &  0 & 0  \\
     0  & u^{-1} & 0 \\
     0 & 0 & 1
\end{pmatrix} $$ \end{enumerate} 

\end{document}